\documentclass[11pt]{article}

\usepackage{latexsym}
\usepackage{amsfonts}
\usepackage{amsmath}
\usepackage{mathrsfs}  
\usepackage{dsfont}    
\usepackage{bbold}     
\usepackage[english]{babel}
\usepackage{caption}
\usepackage{epsfig}
\usepackage{float}
\usepackage{subfigure}

\DeclareMathOperator{\Val}{Val}

\newtheorem{theorem}{Theorem}

\newtheorem{lemma}[theorem]{Lemma}

\newtheorem{hyp}{Assumption}

\newcommand{\zerarcounters}{\setcounter{equation}{0}\setcounter{theorem}{0}}

\newcommand{\ZZZ}{\mathds{Z}}

\newcommand{\RRR}{\mathds{R}}
\newcommand{\TTT}{\mathds{T}}

\newcommand{\MM}{{\mathcal M}}

\newcommand{\calP}{{\mathcal P}}

\newcommand{\TT}{{\mathcal T}}

\newcommand{\gotn}{{\mathfrak n}}

\newcommand{\gotB}{{\mathfrak B}}

\newcommand{\gotN}{{\mathfrak N}}

\newcommand{\gotR}{{\mathfrak R}}
\newcommand{\gotS}{{\mathfrak S}}
\newcommand{\gotT}{{\mathfrak T}}

\newcommand{\Fullbox}{{\rule{2.0mm}{2.0mm}}}

\newcommand{\EP}{\hfill\Fullbox\vspace{0.2cm}}
\newcommand{\prova}{\noindent{\it Proof. }}
\newcommand{\io}{\infty}
\newcommand{\e}{\varepsilon}
\newcommand{\al}{\alpha}

\newcommand{\n}{\nu}

\newcommand{\g}{\gamma}
\newcommand{\om}{\omega}

\newcommand{\s}{\sigma}

\newcommand{\laa}{\langle}
\newcommand{\raa}{\rangle}

\newcommand{\oo}{\boldsymbol{\omega}}
\newcommand{\nn}{\boldsymbol{\nu}}
\newcommand{\pps}{\boldsymbol{\psi}}
\newcommand{\vzero}{\boldsymbol{0}}
\newcommand{\der}{{\rm d}}
\newcommand{\ii}{{\rm i}}

\def\ins#1#2#3{\vbox to0pt{\kern-#2 \hbox{\kern#1 #3}\vss}\nointerlineskip}

\begin{document}

\title{\bf Quasi-periodic motions\\in
strongly dissipative forced systems}

\author
{\bf Guido Gentile
\vspace{2mm}
\\ \small 
$^\dagger$Dipartimento di Matematica, Universit\`a di Roma Tre, Roma,
I-00146, Italy.
\\ \small 
E-mail: gentile@mat.uniroma3.it}

\date{}

\maketitle

\begin{abstract}
We consider a class of ordinary differential equations describing
one-dimensional systems with a quasi-periodic forcing term
and in the presence of large damping. We discuss the conditions
to be assumed on the mechanical force and the forcing term
for the existence of quasi-periodic solutions which have
the same frequency vector as the forcing.
\end{abstract}




\zerarcounters
\section{Introduction}
\label{sec:1}

In this paper we study the same problem considered in \cite{GBD1,GBD2},
that is the existence of quasi-periodic motions in strongly dissipative
forced systems, with the aim of removing as far as possible the
non-degeneracy condition on the mechanical force and the forcing.

We consider one-dimensional systems with a quasi-periodic forcing
term in the presence of strong damping, described
by ordinary differential equations of the form
\begin{equation}
\e \ddot x + \dot x + \e g(x) = \e f(\oo t) ,
\label{eq:1.1} \end{equation}
where $g(x)$ is the \textit{mechanical force},
$f(\oo t)$ is the \textit{forcing term},
$\oo\in\RRR^{d}$ is the \textit{frequency vector} of the forcing,
and $\g=1/\e>0$ is the \textit{damping coefficient}.
Systems of the form (\ref{eq:1.1}) naturally arise in
classical mechanics and electronic engineering;
we refer to \cite{BDGM,GBD1} for physical motivations.
A classical question in the case of forced systems asks for
\textit{response solutions}, that is solutions which are quasi-periodic
with the same frequency vector as the forcing. Note that in
(\ref{eq:1.1}) the forcing is not assumed to be small, as usually
done \cite{St} (see also \cite{BHS} for a review of recent developments):
it is the inverse of the damping coefficient which plays
the role of the perturbation parameter.

Both functions $g$ and $f$ will be assumed to be analytic
in their arguments, with $f$ quasi-periodic, i.e.
\begin{equation}
f(\pps) = \sum_{\nn\in \ZZZ^{d}}
{\rm e}^{\ii\nn \cdot \pps} f_{\nn} , \qquad \pps \in \TTT^{d} ,
\label{eq:1.2} \end{equation}
with average $\laa f \raa = f_{\vzero}$, and
$\cdot$ denoting the scalar product in $\RRR^{d}$.
By the analyticity assumption on $f$, one has
$|f_{\nn}| \le \Phi {\rm e}^{-\xi|\nn|}$
for suitable positive constants $\Phi$ and $\xi$.

A Diophantine condition is assumed on $\oo$. Define
the \textit{Bryuno function} \cite{B}
\begin{equation}
\gotB(\oo) = \sum_{n=0}^{\infty} \frac{1}{2^{n}} \log
\frac{1}{\al_{n}(\oo)} , \qquad \al_{n}(\oo) = \inf
\{ |\oo\cdot\nn| : \nn \in \ZZZ^{d}
\hbox{ such that } 0<|\nn|\le 2^{n} \} .
\label{eq:1.3} \end{equation}
%
\begin{hyp}\label{hyp:1} 
The frequency vector $\oo$ satisfies
the Bryuno condition $\gotB(\oo)<\infty$.
\end{hyp}

Note that if $\oo\in\RRR^{d}$ satisfies the
standard Diophantine condition
$\left|\oo\cdot\nn \right| \ge \g_{0} |\nn|^{-\tau}$
for all $\nn\in\ZZZ^{d}_{*}$,
where $|\nn|:=|\nn|_{1}\equiv|\n_{1}|+\ldots+|\n_{d}|$ and
$\ZZZ^{d}_{*} := \ZZZ^{d} \setminus\{\vzero\}$, and
for some positive constants $\g_{0}$ and $\tau$, then
it also satisfies (\ref{eq:1.3}).
Recently, the Bryuno condition has received a lot of attention
in the theory of small divisor problems; see for instance
\cite{MMY,G2,G3,L,P} and papers cited therein.

The following assumption will be made on the functions $g$ and $f$.

\begin{hyp}\label{hyp:2} 
There exists $c_{0}\in\RRR$ such that $x=c_{0}$ is a zero
of odd order $\gotn$ of the equation
\begin{equation}
g(x) - f_{\vzero} = 0 ,
\label{eq:1.4} \end{equation}
that is $\der^{\gotn}g/\der x^{\gotn}(c_{0})\neq0$ and, if $\gotn>1$,
$\der^{k} g/\der x^{k}(c_{0})=0$ for $k=1,\ldots,\gotn -1$.
\end{hyp}

Of course, for given force $g(x)$, one can read Assumption \ref{hyp:2}
as a condition on the forcing term.

In \cite{GBD1,GBD2} we considered Assumption \ref{hyp:2}
with $\gotn=1$, and, in that case, we proved that
for $\e>0$ small enough there exists a quasi-periodic solution
with frequency vector $\oo$, reducing to $c_{0}$ as $\e$ tends to $0$,
and that such a solution is analytic in a circle tangent
at the origin to the vertical axis.

In this paper we show that the same result of existence extends
under the weaker Assumption \ref{hyp:2}. We also show
that in the case of even $\gotn$ a quasi-periodic solution
oscillating around $x=c_{0}$ fails to exist.
More formal statements are given in Section \ref{sec:2}.

The paper is organised as follows. In Section \ref{sec:2}
we split the problem into two equations, which, using
standard terminology, will be called the \textit{range equation}
and the \textit{bifurcation equation}. The first one involves
small denominator problems, and will be solved iteratively
in Section \ref{sec:3} by using techniques of multiscale analysis
\cite{G1,G2,G3}; from a technical point of view this is the core
of the paper. The second one is an implicit function equation,
and will be discussed in Section \ref{sec:4}.
In Section \ref{sec:5} we show that in the case of zeroes
of even order for the equation (\ref{eq:1.4}), a quasi-periodic
solution of the form $x(t)=c_{0}+O(\e)$ does not exist.
Finally in Section \ref{sec:6} we draw some conclusions and remarks.
The paper is fully self-contained,
and no acquaintance with previous works is required.

\zerarcounters
\section{Setting the problem}
\label{sec:2}

We are interested in the existence of a quasi-periodic solution
with frequency vector $\oo$, hence we expand $x$ as
\begin{equation}
x(t) = c + X(\oo t) , \qquad
X(\pps) = \sum_{\nn\in\ZZZ^{d}_{*}} {\rm e}^{\ii\nn\cdot\pps} X_{\nu} ,
\label{eq:2.1} \end{equation}
where $c=x_{\vzero}$ is the average of $x$ (hence $X$ is
a zero-average function).
Thus, we can rewrite (\ref{eq:1.1}) in Fourier space as
\begin{equation}
\begin{cases}
\left( \ii\oo\cdot\nn \right)
\left( 1 + \ii\e\oo\cdot\nn \right) X_{\nn} +
\e [g \circ (c + X)]_{\nn} = \e f_{\nn} , & \nn \neq \vzero , \\
[g \circ (c + X)]_{\vzero} = f_{\vzero} , & \nn =\vzero ,
\end{cases}
\label{eq:2.2} \end{equation}
where $[F]_{\nn}$ denotes the $\nn$-th Fourier coefficient
of the function $F$.

We shall adopt the following strategy. We shall look for
a solution $c+X$ of the range equation, i.e. the first equation
in (\ref{eq:2.2}), with $c$ arbitrary, and thereafter we shall fix $c$
in such a way that the bifurcation equation, i.e. the second equation
in (\ref{eq:2.2}), be satisfied. This suggests us to consider,
besides the equations (\ref{eq:2.2}), also the equation
\begin{equation}
\left( \ii\oo\cdot\nn \right)
\left( 1 + \ii\e\oo\cdot\nn \right) X_{\nn} +
\e [g \circ (c + X)]_{\nn} = \e f_{\nn} , \qquad \nn \neq \vzero ,
\label{eq:2.3} \end{equation}
and we shall look for a solution of the form (\ref{eq:2.1}) to
(\ref{eq:2.3}). In Section \ref{sec:3} we shall prove that, for
any $c\in\RRR$ close enough to $c_{0}$ and all $\e$ small enough,
there exist such a solution $x=c+X(\oo t;\e,c)$.
Then in Section \ref{sec:4} we shall study the bifurcation equation
\begin{equation}
[g(c+X(\cdot;\e,c)]_{\vzero} = f_{\vzero} ,
\label{eq:2.4} \end{equation}
and we shall see that for $\e$ small enough there exists a solution $c$
to (\ref{eq:2.4}), tending to $c_{0}$ as $\e$ tends to $0$.
More precisely we shall prove the following result.

\begin{theorem} \label{thm:1}
Under the Assumptions \ref{hyp:1} and \ref{hyp:2} for the
ordinary differential equation (\ref{eq:1.1}), for all
$\e$ small enough there exist a continuous function $c(\e)$ and
a response solution $x(t)=c(\e)+X(\oo t;\e,c(\e))$ to (\ref{eq:1.1}),
with $c(0)=c_{0}$ and the function $X(\pps;\eta,c)$ which is
$C^{\io}$ in $\eta$ and $c$, vanishing at $\eta=0$,
and $2\pi$-periodic, analytic and zero-average in $\pps$.
\end{theorem}

If $c_{0}$ is not a zero of the equation (\ref{eq:1.4}), obviously
there is no quasi-periodic solution to (\ref{eq:1.1}) reducing to
$c_{0}$ as $\e$ tends to $0$. We shall show that the same non-existence
result holds if $c_{0}$ is a zero of even order of (\ref{eq:1.4}).
Therefore, the following result strengthens Theorem \ref{thm:1}.

\begin{theorem} \label{thm:2}
Under Assumption \ref{hyp:1} for the ordinary differential equation
(\ref{eq:1.1}), there exists a quasi-periodic solution
of the form as in Theorem \ref{thm:1} if and only if
Assumption \ref{hyp:2} is satisfied. In particular if $c_{0}$ is a
zero of even order of (\ref{eq:1.4}) such a solution does not exists.
\end{theorem}

\zerarcounters
\section{The small denominator equation and multiscale analysis}
\label{sec:3}

A graph is a connected set of points and lines.
A \textit{tree} $\theta$ is a graph with no cycle,
such that all the lines are oriented toward a unique
point (\textit{root}) which has only one incident line (root line).
All the points in a tree except the root are called \textit{nodes}.
The orientation of the lines in a tree induces a partial ordering 
relation ($\preceq$) between the nodes. Given two nodes $v$ and $w$,
we shall write $w \prec v$ every time $v$ is along the path
(of lines) which connects $w$ to the root.

We call $E(\theta)$ the set of \textit{end nodes} in $\theta$,
that is the nodes which have no entering line, and $V(\theta)$
the set of \textit{internal nodes} in $\theta$, that is the set of
nodes which have at least one entering line. Set $N(\theta)=
E(\theta) \amalg V(\theta)$. With each end node $v$ we associate
a \textit{mode} label $\nn_{v}\in\ZZZ^{d}_{*}$. For all $v\in N(\theta)$
denote with $s_{v}$ the number of lines entering the node $v$.

We denote with $L(\theta)$ the set of lines in $\theta$. Since a
line $\ell$ is uniquely identified with the node $v$ which it leaves,
we may write $\ell = \ell_{v}$. With each line $\ell$ we associate
a \textit{momentum} label $\nn_{\ell} \in \ZZZ^{d}_{*}$ and
a \textit{scale} label $n_{\ell}\in\ZZZ_{+}$.

The modes of the end nodes and and the momenta of the lines
are related as follows: if $\ell = \ell_{v}$ one has
\begin{equation}
\nn_{\ell} = \sum_{w \in E(\theta) : w \preceq v} \nn_{w} .
\label{eq:3.1} \end{equation}
If $v$ is an internal node then (\ref{eq:3.1}) gives
$\nn_{\ell}=\nn_{\ell_{1}}+\ldots+\nn_{\ell_{s_{v}}}$, where
$\ell_{1},\ldots,\ell_{s_{v}}$ are the lines entering $v$.

We call \textit{equivalent} two trees which can be transformed into
each other by continuously deforming the lines in such a way that
they do not cross each other. Let $\TT_{k,\nn}$ be the set of
inequivalent trees of order $k$ and total momentum $\nn$,
that is the set of inequivalent trees $\theta$ such that $|N(\theta)|=
|V(\theta)|+|E(\theta)|=k$ and the momentum of the root line is $\nn$.

A cluster $T$ on scale $n$ is a maximal set of nodes and lines
connecting them such that all the lines have scales $n'\le n$
and there is at least one line with scale $n$. The lines entering
the cluster $T$ and the possible line coming out from it (unique if
existing at all) are called the external lines of the cluster $T$.
Given a cluster $T$ on scale $n$, we shall denote by $n_{T}=n$ the
scale of the cluster. We call $V(T)$, $E(T)$, and $L(T)$ the set of
internal nodes, of end nodes, and of lines of $T$, respectively;
note that the external lines of $T$ do not belong to $L(T)$.

We call self-energy cluster any cluster $T$ such that
$T$ has only one entering line $\ell_{T}^{2}$ and one exiting
line $\ell_{T}^{1}$, and one has $\sum_{v\in E(T)} \nn_{v} = \vzero$
(and hence $\nn_{\ell_{T}^{1}}=\nn_{\ell_{T}^{2}}$).
Call $\calP_{T}$ the path of lines $\ell\in L(T)$ connecting
$\ell_{T}^{2}$ to $\ell_{T}^{1}$, and set $x_{T}=
\oo\cdot\nn_{\ell_{T}^{1}}=\oo\cdot\nn_{\ell_{T}^{2}}$.
Let $\gotT_{k,\nn}$ be the set of \textit{renormalised trees} in
$\TT_{k,\nn}$, i.e. of trees in $\TT_{k,\nn}$ which do not
contain any self-energy clusters.

If we write
\begin{equation}
g(x) = \sum_{s=0}^{\io} g_{s}(c) (x-c)^{s} ,
\qquad g_{s}(c) := \frac{1}{s!} \frac{\der^{s}}{\der x^{s}} g(c) ,
\label{eq:3.2} \end{equation}
then we can choose $r>0$ such that $|g_{s}(c)| \le \Gamma^{s}$ for
all $c\in B_{r}(c_{0})$, with the constant $\Gamma$ independent of $c$.

Let $\psi$ be a non-decreasing $C^{\infty}$ function defined
in $\RRR_{+}$, such that
\begin{equation}
\psi(u) = \left\{
\begin{array}{ll}
1 , & \text{for } u \geq 1 , \\
0 , & \text{for } u \leq 1/2 ,
\end{array} \right.
\label{eq:3.3} \end{equation}
and set $\chi(u) := 1-\psi(u)$. For all $n \in \ZZZ_{+}$ define
$\chi_{n}(u) := \chi(u/4\al_{n}(\oo))$ and $\psi_{n}(u) :=
\psi(u/4\al_{n}(\oo))$, and set
\begin{equation}
\Xi_{n}(x)=\chi_{0}(|x|)\ldots \chi_{n-1}(|x|) \chi_{n}(|x|) , \qquad
\Psi_{n}(x)=\chi_{0}(|x|)\ldots \chi_{n-1}(|x|) \psi_{n}(|x|) .
\label{eq:3.4} \end{equation}

We associate with each node $v$ a \textit{node factor}
\begin{equation}
F_{v} = \begin{cases}
- \displaystyle{ \frac{1}{s_{v}!} g_{s_{v}}(c) } ,
& v \in V(\theta) , \\
f_{\nn_{v}} , & v \in E(\theta) , \end{cases}
\label{eq:3.5} \end{equation}
and we associate with each line $\ell$ a \textit{propagator}
\begin{equation}
G_{\ell} = G^{[n_{\ell}]}(\oo\cdot\nn_{\ell};\e,c) ,
\label{eq:3.6} \end{equation}
where the functions $G^{[n]}(x;\e,c)$
are recursively defined for $n\ge 0$ as
\begin{subequations}
\begin{align}
& \hskip-.3truecm
G^{[n]}(x;\e,c) =
\frac{\Psi_{n}(x)}{ \ii x(1+\ii \e x) - \MM^{[n-1]}(x;\e,c)} ,
\label{eq:3.7a} \\
& \hskip-.3truecm
\MM^{[n]}(x;\e,c) = \MM^{[n-1]}(x;\e,c) +
\Xi_{n}(x)
M^{[n]}(x;\e,c) , \quad
M^{[n]}(x;\e,c) = \!\! \sum_{T\in\gotR_{n}}\Val(T,x;\e,c) ,
\label{eq:3.7b}
\end{align}
\label{eq:3.7} \end{subequations}
\vskip-.3truecm
\noindent where $\MM^{[-1]}(x;\e,c)=\e g_{1}(c)$, $\gotR_{n}$ is the set
of renormalised self-energy clusters, i.e. of self-energy clusters which
do not contain any further self-energy clusters, on scale $n$, and
\begin{equation}
\Val(T,x;\e,c) = \Big(\prod_{\ell \in L(T)} G_{\ell} \Big) 
\Big( \prod_{v \in N(T)} F_{v} \Big)
\label{eq:3.8} \end{equation}
is called the value of the self-energy cluster $T$.
Note that $\MM^{[-1]}(x;\e,c)=0$ for $\gotn>1$ in Assumption \ref{hyp:2}.

Set
\begin{equation}
X_{\nn}^{[k]} = \sum_{\theta \in \gotT_{k,\n}}
\Val(\theta;\e,c) , \qquad
\Val(\theta;\e,c) = \Big(\prod_{\ell \in L(\theta)} G_{\ell} \Big) 
\Big( \prod_{v \in N(\theta)} F_{v} \Big) ,
\label{eq:3.9} \end{equation}
where $\Val(\theta;\e,c)$ is called the value of the tree $\theta$,
and define the \textit{renormalised  series}
\begin{equation}
\overline X(\pps;\e,c)= \sum_{\nn\in\ZZZ^{d}_{*}}
{\rm e}^{\ii\nn\cdot\pps} \overline X_{\nn} ,
\qquad \overline X_{\nn} = \sum_{k=1}^{\io} \e^{k} X_{\nn}^{[k]} .
\label{eq:3.10} \end{equation}
Set also
\begin{equation}
M(\theta) = \sum_{v \in E(\theta)} |\nn_{v}| , \qquad
M(T) = \sum_{v \in E(T)} |\nn_{v}| , \qquad
\label{eq:3.11} \end{equation}
and call $\gotN_{n}(\theta)$ the number of lines $\ell\in L(\theta)$
such that $n_{\ell}\ge n$, and $\gotN_{n}(T)$ the number of lines
$\ell\in L(T)$ such that $n_{\ell}\ge n$.

Finally define
\begin{equation}
n(\nn) = \inf \left\{ n \in \ZZZ_{+} : |\nn| \le 2^{n} \right\} .
\label{eq:3.12} \end{equation}
Note that $|\oo\cdot\nn| \ge \al_{n(\nn)}(\oo)$, and 
$\al_{n'}(\oo)<\al_{n}(\oo)$ implies $n'>n$.

\begin{lemma} \label{lem:1}
For any renormalised tree $\theta$,
one has $\gotN_{n}(\theta) \le 2^{-(n-2)}M(\theta)$.
\end{lemma}

\prova We prove that $\gotN_{n}(\theta) \le \max \{ 0, 2^{-(n-2)}
M(\theta)-1\}$ by induction on the number of nodes of $\theta$.
If $N(\theta)=1$ and $\gotN_{n}(\theta)=1$, then $\theta$
has only one line $\ell$ and $n_{\ell}\ge n$. Thus,
$|\oo\cdot\nn_{\ell}| \le \al_{n-1}(\oo)/4$, so that
$n(\nn_{\ell}) \ge n$, and hence $|\nn_{\ell}| > 2^{n-1}$,
which implies $2^{-(n-2)}M(\theta)=2^{-(n-2)}|\nn_{\ell}| \ge 2$.

If $N(\theta) > 1$, let $\ell_{0}$ be the root line of $\theta$
and set $\nn=\nn_{\ell_{0}}$. If $n_{\ell_{0}} < n$ the assertion follows
from the inductive hypothesis.
If $n_{\ell_{0}}\ge n$, call $\ell_{1},\ldots,\ell_{m}$ the lines
with scale $\ge n$ which are closest to $\ell_{0}$. The case $m=0$
is trivial. If $m\ge 2$ the bound follows once more from the
inductive hypothesis. Finally, if $m=1$, then $\ell_{1}$ is
the entering line of a cluster $T$ and $\nn'\neq \nn$,
where $\nn'= \nn_{\ell_{1}}$.
Then $|\oo\cdot(\nn-\nn')| \le \al_{n-1}(\oo)/2$, so that
$n(\nn-\nn')\ge n-1$, and hence $M(T)\ge |\nn-\nn'| > 2^{n-2}$.
Therefore, if $\theta_{1}$ is the tree with root line $\ell_{1}$,
one has $M(\theta)=M(T)+M(\theta_{1})$ and hence
\begin{equation}
\gotN_{n}(\theta) = 1 + \gotN_{n}(\theta_{1}) \le
2^{-(n-2)}M(\theta_{1}) \le
2^{-(n-2)}M(\theta) - 2^{-(n-2)}M(T) \le 2^{-(n-2)}M(\theta) - 1 .
\nonumber \end{equation}
Therefore the assertion follows also in this case.\EP

\begin{lemma} \label{lem:2}
Assume there exists a constant $C_{0}$ such that
$|G^{[n]}(x;\e,c)|\le C_{0}/\al_{n}(\oo)$ for all $n\in\ZZZ_{+}$.
Then there exists $\e_{0}>0$ such that,
for all $c\in B_{r}(c_{0})$ and all $|\e|<\e_{0}$,
the series $c+\overline X(\oo t;\e,c)$ converges.
\end{lemma}

\prova Set $D_{0}=\max\{\Gamma,\Phi\}$.
By assumption for all $\theta\in\gotT_{\nn,k}$ one has
\begin{eqnarray}
\left| \Val(\theta;\e,c) \right|
& \!\!\! \le \!\!\! &
C_{0}^{k} D_{0}^{k} {\rm e}^{-\xi M(\theta)} 
\Big(\prod_{\ell \in L(\theta)} \al_{n_{\ell}}^{-1}(\oo) \Big) 
\le C_{0}^{k} D_{0}^{k} {\rm e}^{-\xi M(\theta)} \al_{n_{0}}^{-k}(\oo)
\prod_{n=n_{0}+1}^{\io} {\rm e}^{\gotN_{n}(\theta) \log 1/\al_{n}(\oo)}
\nonumber
\\
& \!\!\! \le \!\!\! &
C_{0}^{k} D_{0}^{k} {\rm e}^{-\xi M(\theta)} \al_{n_{0}}^{-k}(\oo)
\exp \left( 4 M(\theta) \sum_{n=n_{0}+1}^{\io} \frac{1}{2^{n}}
\log \frac{1}{\al_{n}(\oo)} \right) ,
\nonumber \end{eqnarray}
for arbitrary $n_{0}\in\ZZZ_{+}$. The last sum converges
by Assumption \ref{hyp:1}, so that one can choose $n_{0}$ such that
\begin{equation}
\left| \Val(\theta;\e,c) \right| \le 
C_{0}^{k} D_{0}^{k} \al_{n_{0}}^{-k}(\oo) {\rm e}^{-\xi' M(\theta)} ,
\nonumber \end{equation}
with $\xi'=\xi/2$. This is enough to prove the lemma.\EP

\begin{lemma} \label{lem:3}
For any self-energy cluster $T\in\gotR_{n}$ such that
$\Xi_{n}(x_{T}) \neq 0$, one has
$M(T) \ge 2^{n-1}$ and $\gotN_{p}(T) \le 2^{-(p-2)}M(T)$
for all $p \le n$.
\end{lemma}

\prova We first prove the bound $M(T) \ge 2^{n-1}$ for $T\in\gotR_{n}$
such that $\Xi_{n}(x_{T}) \neq 0$.
By construction any $T\in\gotR_{n}$ has at least one line $\ell$
with scale $n_{\ell}=n$. If $\ell\notin\calP_{T}$ then $\ell$ is the root
line of a tree $\theta$ such that $\gotN_{n}(\theta)\le
2^{-(n-2)}M(\theta)$ by Lemma \ref{lem:1}, so that
$1 \le \gotN_{n}(\theta) \le 2^{-(n-2)}M(T)$, which yields the bound. 
If all lines with scale $n$ are along $\calP_{T}$ then call
$\ell$ that which is closest to $\ell_{T}^{2}$: by construction
$\ell_{T}^{2}$ and $\ell$ are the entering line and the exiting line,
respectively, of a cluster $T'\subset T$, and
$|\nn_{\ell}-\nn_{\ell_{T}^{2}}| \le M(T')$.
Moreover one has $|\oo\cdot\nn_{\ell}|,|\oo\cdot\nn_{\ell_{T}^{2}}|
\le \al_{n-1}(\oo)/4$, hence $M(T) \ge M(T') \ge
|\nn_{\ell}-\nn_{\ell_{T}^{2}}|\ge 2^{n-1}$.

Now we prove that for $T\in\gotR_{n}$ such that $\Xi_{n}(x_{T})\neq0$
one has $\gotN_{p}(T) \le 2^{-(p-2)}M(T)$ for all $p\le n$.
More generally we prove the bound for the elements
of a wider class of graphs. We say that a subset 
$\widetilde T$ of a tree belongs to the class $\gotS_{n,p}$ if
$\widetilde T$ has one exiting line $\ell_{\widetilde T}^{1}$
and one entering line $\ell_{\widetilde T}^{2}$,
both with scale $\ge p$, and all lines $\ell$ in $\widetilde T'$
have scale $n_{\ell} \le n$. Then we prove the bound
$\gotN_{p}(\widetilde T) \le 2^{-(p-2)}M(\widetilde T)$
for all elements $\widetilde T$ of the class $\gotS_{n,p}$.
The proof is by induction on the number of nodes.
Given a subset $\widetilde T$,
let $\ell_{1},\ldots,\ell_{m}$ the lines on scale $\ge p$ closest
to $\ell_{\widetilde T}^{1}$. If $m=0$ then the bound follows easily.
Also the case in which all lines do not belong to the path
$\calP_{\widetilde T}$ can be easily discussed by relying on
Lemma \ref{lem:1}. If at least one line, say $\ell_{1}$, is along
the path $\calP_{\widetilde T}$, then one has
\begin{equation}
\gotN_{p}(\widetilde T) \le 1 +
\gotN_{p}(\widetilde T') + \gotN_{p}(\theta_{2}) + \ldots +
\gotN_{p}(\theta_{m}) , 
\nonumber \end{equation}
where $\theta_{i}$, $i=2,\ldots,m$, is the tree with root line
$\ell_{i}$, while $\widetilde T'$ is a subset with the same
properties as $\widetilde T$, i.e. inside the the same class
$\gotS_{n,p}$, but with $N(\widetilde T')<N(\widetilde T)$.
Hence, by the inductive hypothesis,
one has $\gotN_{p}(\widetilde T') \le 2^{-(p-2)}M(\widetilde T')$.
Then the assertion follows once more.
To conclude the proof simply note that if $T\in\gotR_{n}$ then
$T\in\gotS_{n,p}$ for all $p\le n$.\EP

\begin{lemma} \label{lem:4}
Assume the propagators $G^{[p]}(x;\e,c)$ are differentiable in $x$
and there exist constants $C_{0}$ and $C_{1}$ such that
$|G^{[p]}(x;\e,c)| \le C_{0}/\al_{p}(\oo)$ and
$|\partial_{x} G^{[p]}(x;\e,c)| \le C_{1}/\al_{p}^{3}(\oo)$ for all $p<n$.
Then there exists $\e_{0}>0$ such that, for all $c\in B_{r}(c_{0})$
and all $|\e|<\e_{0}$, the function $x \mapsto M^{[n]}(x;\e,c)$
is differentiable, and one has
\begin{equation}
\left| M^{[n]} (x;\e,c) \right| ,
\left| \partial_{x} M^{[n]} (x;\e,c) \right| \le
D_{1} |\e|^{2} {\rm e}^{-D_{2}2^{n}} ,
\nonumber \end{equation}
for some positive constants $D_{1}$ and $D_{2}$.
\end{lemma}

\prova By proceeding as in the proof of Lemma \ref{lem:2}, one finds
\begin{equation}
\left| \Val(T,x;\e,c) \right| \le
C_{0}^{k} D_{0}^{k} {\rm e}^{-\xi M(T)} 
\al_{n_{0}}^{-k}(\oo) \exp
\left( 4 M(T) \sum_{n=n_{0}+1}^{\io} \frac{1}{2^{n}}
\log \frac{1}{\al_{n}(\oo)} \right) ,
\nonumber \end{equation}
with $n_{0}$ chosen as in the proof of Lemma \ref{lem:2}.
Then one can use Lemma \ref{lem:3} to bound $M(T)$,
and the observation that any self-energy cluster $T$
has at least two nodes to obtain the factor $\e^{2}$.
This proves the bound on $M^{[n]} (x;\e,c)$.

To obtain the bound on $\partial_{x}M^{[n]} (x;\e,c)$ simply note that
\begin{equation}
\partial_{x} M^{[n]}(x;\e,c) = \sum_{T\in\gotR_{n}}
\Big( \prod_{v \in E(T) \cup V(T)} F_{v} \Big)
\sum_{\ell\in \calP_{T}} \partial_{x} G_{\ell}
\Big(\prod_{\ell' \in L(\theta) \setminus\{\ell\}} G_{\ell'} \Big) ,
\nonumber \end{equation}
where $\partial_{x} G_{\ell}$ can be bounded as $|\partial_{x}
G_{\ell}| \le C_{1}/\al_{n_{\ell}}^{3}(\oo)$ by hypothesis.\EP

\begin{lemma} \label{lem:5}
Assume there exists a constant $C_{0}$ such that
$|G^{[p]}(x;\e,c)|\le C_{0}/\al_{p}(\oo)$ for all $p<n$.
Then one has $(\MM^{[p]}(x;\e,c))^{*} = \MM^{[p]}(-x;\e,c)$
for all $p \le n$.
\end{lemma}

\prova The proof is by induction on $p$. First of all note that
if $(\MM^{[p]}(x;\e,c))^{*} = \MM^{[p]}(-x;\e,c)$ then
$(G^{[p]}(x;\e,c))^{*} = G^{[p]}(-x;\e,c)$ by (\ref{eq:3.7a}).
Moreover one has $F_{v}^{*}=F_{v}$ for all internal nodes $v\in V(T)$ and
$F_{v}^{*}=f_{\nn_{v}}^{*}=f_{-\nn_{v}}$ for all end nodes $v\in E(T)$.

Let $T$ a self-energy cluster contributing to $\MM^{[p]}(x;\e,c)$
-- see (\ref{eq:3.7b}) -- for $p\le n$; then $T\in \gotR_{q}$ for
some $q\le p$. Together with $T\in\gotR_{q}$ consider also the
self-energy cluster $T'\in\gotR_{q}$ obtained from $T$ by changing
the signs of the mode labels of all the end nodes $v\in E(T)$.
Note that there is a one-to-one correspondence between
the self-energy clusters $T$ and $T'$.
The node factors corresponding to the end nodes $v \in E(T')$
become $f_{-\nn_{v}}$, and, if we revert the momentum of the
entering line $\ell_{T'}^{2}$, the momenta of all the lines
$\ell\in L(T')$ also change sign, that is $\nn_{\ell}$ is replaced
with $-\nn_{\ell}$ for all $\ell\in L(T')$.

The definition (\ref{eq:3.8}) and the inductive hypothesis
yield $(\Val(T,x;\e,c))^{*}=\Val(T',-x;\e,c)$ for all $q \le p$
and all $T\in\gotR_{q}$. Then (\ref{eq:3.7b}) implies the assertion.\EP

\begin{lemma} \label{lem:6}
For all $n\in\ZZZ_{+}$ the function $x \mapsto \MM^{[n]}(x;\e,c))$
is differentiable and one has
$|\ii x (1 + \ii \e x) - \MM^{[n]}(x;\e,c)|\ge |x|/2$
for all $c\in B_{r}(c_{0})$ and all $\e$ small enough.
\end{lemma}

\prova The proof is by induction on $n$. Assume that
the functions $x \mapsto \MM^{[p]}(x;\e,c)$ are differentiable and
one has $|\ii x (1 + \ii \e x) - \MM^{[p]}(x;\e,c)|\ge |x|/2$
for all $p<n$. One can easily verify that then also the propagators
$G^{[p]}(x;\e,c)$ are differentiable and satisfy the bounds
$|\partial_{x} G^{[p]}(x;\e,c)| \le C_{1}/\al_{p}^{3}(\oo)$ 
for all $p \le n$ and for some positive constant $C_{1}$. Indeed one has
\begin{equation}
\partial_{x} G^{[p]}(x;\e,c) =
\frac{\partial_{x}\Psi_{p}(x)}{\ii x(1+\ii \e x)-\MM^{[p-1]}(x;\e,c)} -
\frac{\Psi_{p}(x) \left( \ii - 2\e x -\partial_{x}\MM^{[p-1]}
(x;\e,c) \right)}{(\ii x(1+\ii \e x)-\MM^{[p-1]}(x;\e,c))^{2}} ,
\nonumber \end{equation}
where
\begin{eqnarray}
\partial_{x}\Psi_{p}(x)
& \!\!\! = \!\!\! & \sum_{j=0}^{p-1}
\chi_{0}(|x|) \ldots \partial_{x} \chi_{j}(|x|) \ldots \psi_{p}(|x|) +
\chi_{0}(|x|) \ldots \chi_{p-1}(|x|) \partial_{x} \psi_{n}(|x|)
\nonumber \\
& \!\!\! \le \!\!\! &
C \sum_{j=0}^{p} \al_{j}^{-1}(\oo) \le C p \al_{p}^{-1}(\oo)
\nonumber \end{eqnarray}
for some constant $C$, and
\begin{eqnarray}
\partial_{x}\MM^{[p-1]}(x;\e,c)
& \!\!\! = \!\!\! & \sum_{j=0}^{p-1} \left( \sum_{i=0}^{j} 
\chi_{0}(|x|) \ldots \partial_{x} \chi_{i}(|x|) \ldots \chi_{j}(|x|)
M^{[j]}(x;\e,c) + \Xi_{j}(x) \partial_{x} M^{[j]}(x;\e,c) \right)
\nonumber \\
& \!\!\! \le \!\!\! &
C |\e|^{2}  \sum_{j=0}^{p-1} {\rm e}^{-D_{2}2^{j}}
\left( j \al_{j}^{-1}(\oo) + 1 \right) \le
C' |\e|^{2} p^{2} \al_{p}^{-1}(\oo) ,
\nonumber \end{eqnarray}
for some constants $C,C'$.

Then we can apply Lemma \ref{lem:4} to conclude that $\MM^{[n]}(x;\e,c)$
is differentiable and its derivative with respect to $x$
is accordingly bounded. Therefore
\begin{equation}
\ii x (1 + \ii \e x) - \MM^{[n]}(x;\e,c) =
\ii x (1 + \ii \e x) - \MM^{[n]}(0;\e,c) -
\left( \MM^{[n]}(x;\e,c) - \MM^{[n]}(0;\e,c) \right) ,
\nonumber \end{equation}
where $\MM^{[n]}(0;\e,c)$ is real by Lemma \ref{lem:6}, and
\begin{equation}
\left| \MM^{[n]}(x;\e,c)) - \MM^{[n]}(0;\e,c) \right|
\le C |\e|^{2} |x| ,
\nonumber \end{equation}
for some constant $C$, by Lemma \ref{lem:4}.\EP

\begin{lemma} \label{lem:7}
Then there exists $\e_{0}>0$ such that,
for all $c\in B_{r}(c_{0})$ and all $|\e|<\e_{0}$, the function
$c+\overline X(\oo t;\e,c)$ solves (\ref{eq:2.3}).
\end{lemma}

\prova We have to prove that the coefficients $\overline X_{\nn}$,
defined abstractly through (\ref{eq:3.10}), solve the
first equation in (\ref{eq:2.2}), i.e.
\begin{equation}
\left( \ii\oo\cdot\nn \right)
\left( 1 + \ii\e\oo\cdot\nn \right) \overline X_{\nn} +
\e [g \circ (c + \overline X)]_{\nn} = \e f_{\nn} ,
\qquad \nn \neq \vzero .
\nonumber \end{equation}
Set $D_{n}(x;\e,c)= \ii x \left( 1+\ii\e x \right)-\MM^{[n]}(x;\e,c)$,
so that $G^{[n]}(x;\e,c)= \Psi_{n}(x)/D_{n}(x;\e,c)$, and
$G(x)=1/(\ii x)(1+\ii\e x)$. Write also
\begin{equation}
\overline X_{\nn} = \sum_{n=0}^{\io} \overline X_{\nn,n} ,
\qquad \overline X_{\nn,n} =
\sum_{k=1}^{\io} \e^{k} \sum_{\theta\in\gotT_{k,\nn,n}}
\Val(\theta;\e,c) ,
\nonumber \end{equation}
where $\gotT_{k,\nn,n}$ is the subset of $\gotT_{k,\nn}$
of the renormalised trees with root line with scale $n$.

If we define
\begin{equation}
\Omega(\nn,\e,c) =
G(\oo\cdot\nn) \left[ \e f - \e g(c +
\overline X(\cdot;\e,c) \right]_{\nn} ,
\label{eq:3.13} \end{equation}
then we have to prove that $\Omega(\nn,\e,c) = \overline X_{\nn}$.

By setting
\begin{equation}
\Psi_{j,n}(x) = \chi_{j}(|x|)\ldots \chi_{n-1}(|x|) \psi_{n}(|x|) ,
\quad n>j , \qquad \qquad \Psi_{n,n}(x) = \psi_{n}(|x|) , 
\label{eq:3.14} \end{equation}
note that
\begin{equation}
\Psi_{0,n}(x) = \Psi_{n}(x), \qquad \qquad
\sum_{n=j}^{\io} \Psi_{j,n}(x) = 1 \quad \forall j \ge 0 .
\label{eq:3.15} \end{equation}
Then, by using the last identity in (\ref{eq:3.15}) with $j=0$,
we can rewrite (\ref{eq:3.13}) as
\begin{equation}
\Omega(\nn,\e,c) =
G(\oo\cdot\nn) \sum_{n=0}^{\io}
D_{n}(\oo\cdot\nn;\e,c) \, G^{[n]}(\oo\cdot\nn;\e,c)
\left[ \e f - \e g(c + \overline X(\cdot;\e,c)) \right]_{\nn} ,
\label{eq:3.16} \end{equation}
where we can expand
\begin{eqnarray}
& & G^{[n]}(\oo\cdot\nn;\e,c)
\left[ \e f_{\nn} - \e g(c + \overline X(\cdot;\e,c)) \right]_{\nn} =
\sum_{k=1}^{\io} \e^{k} \sum_{\theta\in\gotT_{k,\nn,n}} \Val(\theta;\e,c)
\nonumber \\
& & \qquad \qquad
+ \; G^{[n]}(\oo\cdot\nn;\e,c)
\sum_{p=n}^{\io} \sum_{j=0}^{n-1} M^{[j]}(\oo\cdot\nn;\e,c)
\sum_{k=1}^{\io} \e^{k} \sum_{\theta\in\gotT_{k,\nn,p}} \Val(\theta;\e,c)
\nonumber \\
& & \qquad \qquad
+ \; G^{[n]}(\oo\cdot\nn;\e,c)
\sum_{p=0}^{n-1} \sum_{j=0}^{p-1} M^{[j]}(\oo\cdot\nn;\e,c)
\sum_{k=1}^{\io} \e^{k} \sum_{\theta\in\gotT_{k,\nn,p}} \Val(\theta;\e,c) ,
\nonumber \end{eqnarray}
where the sum in the second line is present only if $n\ge 1$
and the sum in the third line is present only if $n\ge 2$.
Therefore we obtain
\begin{eqnarray}
& & \Omega(\nn,\e,c) =
G(\oo\cdot\nn) \sum_{n=0}^{\io} D_{n}(\oo\cdot\nn;\e,c)
\overline X_{\nn,n}
\nonumber \\
& & \qquad \qquad
+ \; G(\oo\cdot\nn) \sum_{n=1}^{\io} \Psi_{n}(x)
\sum_{p=n}^{\io} \sum_{j=0}^{n-1} M^{[j]}(\oo\cdot\nn;\e,c)
\overline X_{\nn,p}
\nonumber \\
& & \qquad \qquad
+ \; G(\oo\cdot\nn) \sum_{n=2}^{\io} \Psi_{n}(x)
\sum_{p=0}^{n-1} \sum_{j=0}^{p-1} M^{[j]}(\oo\cdot\nn;\e,c)
\overline X_{\nn,p} .
\nonumber \end{eqnarray}
The second and third lines, summed together, give
\begin{equation}
G(\oo\cdot\nn) \sum_{n=1}^{\io} \overline X_{\nn,n}
\sum_{j=0}^{n-1} M^{[j]}(\oo\cdot\nn;\e,c)
\sum_{p=j+1}^{\io} \Psi_{p}(x) , \qquad \hbox{where} \qquad
\sum_{p=j+1}^{\io} \Psi_{p}(x) = \Xi_{j}(x) ,
\nonumber \end{equation}
where we have written $\Psi_{p}=\Xi_{j}\Psi_{j+1,p}$
and used (\ref{eq:3.14}) to obtain the last equality, so that
(\ref{eq:3.16}) gives
\begin{equation}
\Omega(\nn,\e,c) =
G(\oo\cdot\nn) \sum_{n=0}^{\io}
\left( D_{n}(\oo\cdot\nn;\e,c) + \MM^{[n]}(\oo\cdot\nn;\e,c) 
\right) \overline X_{\nn,n} = \sum_{n=0}^{\io}
\overline X_{\nn,n} = \overline X_{\nn} ,
\nonumber \end{equation}
which proves the assertion.\EP

\begin{lemma} \label{lem:8}
The function $\overline X(\pps;\e,c)$ is $C^{\io}$
in $\e$ and $c$ for $\e$ and $c-c_{0}$ small enough.
\end{lemma}

\prova The previous results imply that $\overline X(\cdot;\e,\cdot)$
is a well defined function of $\e$ for $\e$ small enough. 
By looking at the tree expansion (\ref{eq:3.9}) for the coefficients
$\overline X^{[k]}_{\nn}$ of  $\overline X(\pps;\e,c)$,
one sees that the function depends on $\e$ through the
factors $\e^{k}$ in (\ref{eq:3.10}) and through the propagators
$G_{\ell}$. The first dependence is trivial, and poses no obstacle
in differentiating. Also the dependence through the propagators
can be easily handled thanks to Lemma \ref{lem:6}, which allows
to bound from below the denominators. In particular
for all $m \ge 0$ one finds
\begin{equation}
\left| \partial_{\e}^{m} G^{[p]}(x;\e,c) \right| \le K_{m} /\al_{p}^{m}(\om)
\nonumber \end{equation}
for suitable constants $K_{m}$. Smoothness in $c$ can be discussed
in a similar way, by using analyticity of the force $g$ and again
Lemma \ref{lem:6}. \EP

\zerarcounters
\section{The implicit function equation}
\label{sec:4}

We are left with the implicit function equation (\ref{eq:2.4}), which
can be trivially solved under Assumption \ref{hyp:2}. If we define
\begin{equation}
\Gamma(\e,c) = [g(c+X(\cdot;\e,c)]_{\vzero} - f_{\vzero} ,
\label{eq:4.1} \end{equation}
then the following result holds.

\begin{lemma} \label{lem:9}
There exists a neighbourhood $U \times V$ of $(\e,c)=(0,c_{0})$
such that for all $\e\in U$ there is at least one value $c=c(\e) \in V$,
depending continuously on $\e$, for which one has $\Gamma(\e,c(\e))=0$.
\end{lemma}

\prova Since $\Gamma(0,c)=g(c)-f_{\vzero}$,
Assumption \ref{hyp:2} implies that
\begin{equation}
\frac{\der^{k}}{\der c^{k}} \Gamma(0,c_{0})=0 \quad
\hbox{for} \quad k=0,1,\ldots,\gotn -1 \quad \hbox{and} \quad
\Gamma_{0} = \frac{\der^{\gotn}}{\der c^{\gotn}} \Gamma(0,c_{0}) \neq 0 .
\label{eq:4.2} \end{equation}
Set $\s_{0}={\rm sign}(\Gamma_{0})$ so that $\s_{0}\Gamma_{0}>0$.
By continuity there are neighbourhoods $U$ and $V=[V_{-},V_{+}]$ of
$\e=0$ and $c=c_{0}$, respectively, such that for all $\e\in U$ one has
$\s_{0}\Gamma(\e,c)>0$ for $c =V_{+}$ and $\s_{0}\Gamma(\e,c)<0$
for $c =V_{-}$. Therefore, there exists a continuous curve $c=c(\e)$
such that $\Gamma(\e,c(\e))=0$.\EP

By collecting together the results of the previous sections and
Lemma \ref{lem:9}, Theorem \ref{thm:1} follows.

\zerarcounters
\section{Zeroes of even order}
\label{sec:5}

In this section we prove the following result, which, together with
Theorem \ref{thm:1}, implies Theorem \ref{thm:2}.

\begin{lemma} \label{lem:10}
Under Assumption \ref{hyp:1} for the ordinary differential equation
(\ref{eq:1.1}), assume also that $c_{0}$ is a zero of even order
of (\ref{eq:1.1}). Then there is no quasi-periodic solution
reducing to $c_{0}$ when $\e$ tends to $0$.
\end{lemma}

\prova The analysis of Section \ref{sec:3} shows that a solution
of the range equation (\ref{eq:2.3}) can be proved to exist
under the only Assumption \ref{hyp:1}. Moreover such a solution
is $C^{\io}$ in both $\e$ and $c$ (cf. Lemma \ref{lem:8}).
Then, we study the bifurcation equation (\ref{eq:2.4}) in the case
$c_{0}$ is a zero of even order of (\ref{eq:1.4}).

If we write $c=c_{0}+\zeta$ and expand the function $g \circ (c+X)$
around $c=c_{0}$, then (\ref{eq:2.4}) gives
\begin{equation}
[g(c+X(\cdot;\e,c)]_{\vzero} - f_{\vzero} =
g_{0} \laa \left( \zeta + X \right)^{\gotn} \raa + 
\laa O\left( \zeta + X \right)^{\gotn+1} \raa = 0 ,
\label{eq:5.1} \end{equation} 
where $\gotn! g_{0}=\der^{\gotn}g/\der x^{\gotn}(c_{0})\neq 0$
and $\laa\cdot\raa$ denotes as usual the Fourier component
with label $\nn=\vzero$. If $\e=O(\zeta)$ then we have
$|g_{0}| \laa \left( \zeta + X \right)^{\gotn} \raa \ge C_{1}\e^{\gotn}$
for some positive constant $C_{1}$, because $\gotn$ is even,
and $O\left( \zeta + X \right)^{\gotn+1} = O(\e^{\gotn+1})$,
so that (\ref{eq:5.1}) cannot be satisfied for $\e$ small enough.

If $\e=o(\zeta)$, then
\begin{equation}
\laa \left( \zeta + X \right)^{\gotn} \raa =
\sum_{k=0}^{\gotn} \left( \begin{matrix} \gotn \\ k \end{matrix} \right)
\zeta^{k} \laa X^{\gotn - k} \raa = \zeta^{\gotn} + o(\zeta^{\gotn})
\label{eq:5.2} \end{equation} 
for $\e$ small enough. On the other hand
$O\left( \zeta + X \right)^{\gotn+1} = O(\zeta^{\gotn+1})$,
and hence once more there is no solution to (\ref{eq:5.1})
because of (\ref{eq:5.2}).
The case $\zeta=o(\e)$ can be discussed in a similar way.\EP

\zerarcounters
\section{Conclusions and open problems}
\label{sec:6}

The analysis of the previous sections shows that under Assumptions
\ref{hyp:1} and \ref{hyp:2} the system described by the ordinary
differential equation (\ref{eq:1.1}) admits a response solution.
Under some mild conditions on $g$ one can prove that
such a solution describes a (local) attractor \cite{BDG}.
It would be interesting to investigate whether the same
result can be obtained by only making Assumption \ref{hyp:2} on $g$
and requiring $\der^{\gotn}g/\der x^{\gotn}(c_{0})>0$.
Even more interesting would be to understand whether the same scenario
persists after removing Assumption \ref{hyp:1} on $\oo$.
The analysis of \cite{BDG} shows that, if there is a quasi-periodic
solution of the form considered in Theorem \ref{thm:1} exists,
then it is an attractor (under some conditions on $g$), but
if $\oo$ does not satisfy any Diophantine condition,
such as the Bryuno condition, then the small divisor problem
can not be handled, and it is very unlikely that the dynamics
can be conjugated to the unperturbed one.

The analysis in Section \ref{sec:5} shows
that, if $c_{0}$ is a zero of even order $\gotn$ for the
equation (\ref{eq:1.4}), then no quasi-periodic
solution of the form considered in Theorem \ref{thm:1} exists.
A natural question in that case is, how the dynamics evolves in time,
and what kind of attractors arise.

Furthermore, Theorem \ref{thm:1} states that for all
$\e$ small enough there is a value $c(\e)$ for the average
of $x(t)$, such that the solution exists, but provides nothing
more than continuity about the dependence of $c(\e)$ on $\e$.
Thus, another question which should deserve further investigation
is, if under some further assumption one can prove some
stronger regularity property for the function $c(\e)$ --
note that analyticity fails to hold even in the case of
periodic forcings \cite{GBD1}. In this direction,
the results of \cite{CG} could provide a possible guideline
(even if in this case the implicit function equation
to be studied is no longer analytic), not only to prove
smoothness but also to provide an algorithm to explicitly
construct the function $c(\e)$. Of course,
under the Assumption \ref{hyp:2} on $\oo$,
independently of the conditions on $g$, we have no hope
to prove Borel summability \cite{S} in $\e$ at the origin.
Indeed, this should require a much stronger Diophantine condition
on $\oo$ \cite{CGGG,GBD2}. 


\end{document}